\documentclass[letterpaper, 10 pt, conference]{ieeeconf}  

\IEEEoverridecommandlockouts

\usepackage{graphicx} 
\usepackage{amsmath} 
\usepackage{amssymb}  

\newtheorem{theorem}{Theorem}
\newtheorem{proposition}{Proposition}
\newtheorem{remark}{Remark}

\newcommand{\R}{\mathbb{R}}

\title{\LARGE \bf
A Less Conservative Sufficient Condition for PID Stabilization of Scalar Second-Order Nonlinear Uncertain Systems}

\author{Senhan Yao
\thanks{
This work has been submitted to the IEEE L-CSS for possible publication.
}
\thanks{This work was supported by the National
Natural Science Foundation of China under Grant No. 12288201. \textit{(Corresponding author: Senhan Yao.)}}
\thanks{Senhan Yao is with the Key Laboratory of Systems and
Control, Academy of Mathematics and Systems Science, Chinese Academy of Sciences, Beijing 100190,
China, and also with the School of Mathematical Sciences, University
of Chinese Academy of Sciences, Beijing 100049, China (e-mail: yaosenhan@amss.ac.cn).}}

\begin{document}

\maketitle
\thispagestyle{empty}
\pagestyle{empty}

\begin{abstract}
This letter studies robust set-point regulation of scalar second-order nonlinear uncertain systems using a classical PID controller with constant gains. The scalar second-order model provides a minimal prototype for nonlinear mechanical and electromechanical dynamics, while its velocity-dependent term captures uncertainties such as physical damping and friction. For a positive velocity-derivative bound, existing Lyapunov sufficient conditions certify fixed-gain PID parameter regions that remain separated from the boundary associated with the necessary condition obtained from the worst-case linear model. To reduce this conservatism, this letter proposes an endpoint-balanced quadratic-plus-integral Lyapunov certificate. The key idea is to choose the quadratic cross-term coefficient so that the mixed-term penalty is balanced at the two endpoints of the admissible effective-damping interval before extracting the scalar PID inequality. The resulting condition guarantees global asymptotic regulation for the full derivative-bounded uncertainty class. When the velocity-derivative bound is positive, the proposed condition certifies a fixed-gain PID region that strictly contains those certified by Zhao--Guo and Zhang--Guo. When this bound is zero, the corresponding boundary coincides with that necessary boundary. At the level of Lyapunov analysis, the construction reduces the uniform mixed-term penalty over the entire effective-damping interval.
\end{abstract}

\section{Introduction}
\label{sec:introduction}

How much uncertainty can feedback control tolerate while still guaranteeing global regulation?
This is a basic question in the capability theory of feedback~\cite{Guo2020Feedback}.
For nonlinear systems, a quantitative formulation is to prescribe derivative-bounded uncertainty classes
and ask whether a fixed feedback architecture stabilizes every admissible system.
The classical PID controller is central to this question
because it uses only the instantaneous, integral, and derivative information of the regulation error and requires no accurate plant model~\cite{AstromHagglund2001,Samad2017}.
A natural benchmark is the scalar second-order uncertain system, which can be viewed as a minimal prototype of Euler--Lagrange-type mechanical dynamics with unknown restoring forces, damping/friction, and structural nonlinearities~\cite{CortezDimarogonas2022,TaoRoyDeSchutterBaldi2024}.
This setting is sufficiently simple to allow explicit PID capability boundaries and has been used in recent studies of nonlinear uncertain systems~\cite{ZhaoGuo2017,ZhangGuo2019,ZhaoGuo2017Capability,ZhaoWangXue2025}.

For the derivative-bounded class $F_{L_1,L_2}$ considered in this letter, the worst-case linear model yields the necessary condition $p>0$, $q>0$, $k>0$, and $pq>k$ in terms of the shifted gains $p=k_p-L_1$, $q=k_d-L_2$, $k=k_i$.
Existing scalar Lyapunov sufficient conditions give explicit lower bounds on $p$ that stay strictly above $p=k/q$ when the velocity-derivative bound is positive.
In particular, the scalar boundaries induced by the sufficient conditions in~\cite{ZhaoGuo2017,ZhangGuo2019} both remain strictly above $p=k/q$ when $L_2>0$; we denote these boundaries by $p_{17}(q)$ and $p_{19}(q)$.
In the Lyapunov framework of this letter, those earlier certificates can be viewed as employing mixed-term bounds that are not balanced over the whole admissible effective-damping interval, leaving room for tightening the certified sufficient boundary.
In a complementary direction, Zhao, Wang, and Xue~\cite{ZhaoWangXue2025} enlarge the class of systems that can be globally regulated by changing the feedback architecture to nonlinear PID-type controllers with state- or error-dependent gains.
In contrast, the present letter keeps the classical constant-gain PID architecture
and focuses on reducing the conservatism of scalar Lyapunov sufficient boundaries within this fixed-gain setting.

This letter proposes an endpoint-balanced quadratic-plus-integral Lyapunov certificate.
After the standard integral term eliminates the dependence of the mixed term on the unknown position-dependent coefficient,
the remaining uncertainty enters only through the effective-damping interval $\psi(y,z)\in[q,q+2L_2]$.
For each auxiliary parameter $\mu \in (0,q)$, we select the quadratic cross-term coefficient so that
the normalized mixed-term penalty is balanced at the two endpoints of the shifted interval $[q-\mu,q+2L_2-\mu]$,
before extracting the scalar PID inequality.
This interval-wise placement yields the new sufficient boundary $p>\Pi(q)$.

The main contributions are as follows.
The proposed condition guarantees global asymptotic regulation for the full derivative-bounded class $F_{L_1,L_2}$.
For $L_2>0$, the newly certified fixed-gain PID region strictly contains the regions obtained from the scalar sufficient conditions in~\cite{ZhaoGuo2017,ZhangGuo2019}.
In this notation, the strict inclusion is expressed as
$\Pi(q)<p_{19}(q)<p_{17}(q)$.
For $L_2=0$, the condition exactly recovers the lower boundary $p=k/q$
associated with this necessary condition.
Thus, for $L_2>0$, the result closes part of the certified gap toward the necessary condition without changing the feedback architecture.
At the level of Lyapunov analysis, the improvement comes from reducing the uniform mixed-term penalty
over the entire effective-damping interval, rather than from selecting a different endpoint
or making the cross coefficient state dependent.

The remainder of this letter is organized as follows.
Section~II formulates the problem.
Section~III presents the main result, namely a new sufficient condition,
and compares the certified boundary with the scalar boundaries of
\cite{ZhaoGuo2017,ZhangGuo2019}.
Section~IV provides the proof.
Section~V interprets the endpoint-balanced Lyapunov mechanism
and discusses its relation to the necessary condition.
Section~VI illustrates the result with a reproducible nonlinear example,
and Section~VII concludes the letter.

\emph{Notation.}
Let $\R$ denote the set of real numbers.
For a $C^1$ scalar function $f$, $f_{x_i}$ denotes its partial derivative with respect to $x_i$.
For a symmetric matrix $P$, $P>0$ means that $P$ is positive definite.
The derivative of a Lyapunov function along closed-loop trajectories is denoted by $\dot V$.
All absolute values are scalar absolute values.
Inequalities between scalar functions or parameter-region boundaries are understood pointwise
over their stated domains, and all infima are taken over the explicitly indicated sets.

\section{Problem Formulation}
\label{sec:problem}

Consider the scalar uncertain second-order system
\begin{equation}
\label{eq:plant}
    \dot x_1=x_2,\qquad
    \dot x_2=f(x_1,x_2)+u,
\end{equation}
where $x_1,x_2,u\in\R$ and $f$ is unknown. The control objective is to regulate $x_1(t)$ to a prescribed set point $y^\ast\in\R$ and drive $x_2(t)$ to zero.

The classical PID controller is written in dynamic form as
\begin{equation}
\label{eq:pid}
\begin{aligned}
\dot \xi(t)&=e(t),\\
u(t)&=k_p e(t)+k_i\xi(t)+k_d\dot e(t),\\
e(t)&=y^\ast-x_1(t),
\end{aligned}
\end{equation}
where $k_p,k_i,k_d$ are constant gains and $\xi$ is the integral state. Let $L_1,L_2\ge 0$ be known constants and define
\begin{equation}
\label{eq:function-class}
\begin{aligned}
F_{L_1,L_2}
:=\{\, f\in C^1(\R^2):\ 
& f_{x_1}(x_1,x_2)\le L_1,\\
& |f_{x_2}(x_1,x_2)|\le L_2,\\
& \forall (x_1,x_2)\in\R^2\,\}.
\end{aligned}
\end{equation}
The derivative bound with respect to $x_1$ is one-sided, as in \cite{ZhaoGuo2017}, whereas the bound with respect to $x_2$ is two-sided.

Assume $k_i>0$ and introduce
\begin{equation}
\label{eq:variables}
\begin{aligned}
x(t)&=\xi(t)+\frac{f(y^\ast,0)}{k_i},\\
y(t)&=e(t),\qquad z(t)=\dot e(t).
\end{aligned}
\end{equation}
Set
\begin{equation}
\label{eq:kpq}
    k:=k_i,\qquad
    p:=k_p-L_1,\qquad
    q:=k_d-L_2,\qquad
    L:=L_2.
\end{equation}
In these variables, the closed-loop system takes the normal form
\begin{equation}
\label{eq:normal}
\begin{aligned}
\dot x&=y,\qquad \dot y=z,\\
\dot z&=-kx-\phi(y)y-\psi(y,z)z.
\end{aligned}
\end{equation}
where
\begin{equation}
\label{eq:bounds}
\phi(y)\ge p,\qquad
q\le \psi(y,z)\le q+2L.
\end{equation}

Indeed, define
\begin{equation}
\label{eq:gdef}
    g(y,z)=f(y^\ast,0)-f(y^\ast-y,-z).
\end{equation}
Then $g(0,0)=0$ and
\begin{equation}
\label{eq:zpre}
    \dot z=g(y,z)-kx-k_p y-k_d z.
\end{equation}
Writing
\begin{equation}
\label{eq:ba}
    g(y,z)=b(y)y+a(y,z)z
\end{equation}
with
\begin{equation}
\label{eq:bdef}
b(y)=
\begin{cases}
g(y,0)/y, & y\ne 0,\\
g_y(0,0), & y=0,
\end{cases}
\end{equation}
and
\begin{equation}
\label{eq:adef}
a(y,z)=
\begin{cases}
\bigl(g(y,z)-g(y,0)\bigr)/z, & z\ne 0,\\
g_z(y,0), & z=0,
\end{cases}
\end{equation}
the mean value theorem and \eqref{eq:function-class} give
\begin{equation}
\label{eq:ba-bounds}
    b(y)\le L_1,\qquad |a(y,z)|\le L_2.
\end{equation}
Thus
\begin{equation}
\label{eq:phipsi}
    \phi(y)=k_p-b(y),\qquad \psi(y,z)=k_d-a(y,z),
\end{equation}
which yields \eqref{eq:bounds}. The apparent singularities in these quotient definitions are removable. Since $g\in C^1(\R^2)$, they can equivalently be written as
\begin{equation}
\label{eq:ba-integral}
    b(y)=\int_0^1 g_y(\theta y,0)\,d\theta,\qquad
    a(y,z)=\int_0^1 g_z(y,\theta z)\,d\theta.
\end{equation}
Thus $b$, $a$, $\phi$, and $\psi$ are continuous. Moreover, the vector field represented by \eqref{eq:normal} is exactly
\[
(y,z,g(y,z)-kx-k_py-k_dz),
\]
which is locally Lipschitz because $g\in C^1(\R^2)$. Classical solutions are therefore unique on their maximal intervals of existence, and the Lyapunov and LaSalle arguments below apply.

\section{Main Result}
\label{sec:main}

The next theorem gives an endpoint-balanced certificate that reduces the scalar sufficient boundary.

\begin{theorem}
\label{thm:main}
Consider the closed-loop normal form \eqref{eq:normal} obtained from \eqref{eq:plant}--\eqref{eq:pid} for some $f\in F_{L_1,L_2}$. Assume
\begin{equation}
\label{eq:basic-positive}
    k>0,\qquad q>0.
\end{equation}
If there exists $\mu\in(0,q)$ satisfying
\begin{equation}
\label{eq:new-condition}
\begin{aligned}
4(\mu p-k)
>
\mu^2
\Big(
&\sqrt{q+2L-\mu}-\sqrt{q-\mu}
\Big)^2,
\end{aligned}
\end{equation}
then the origin of \eqref{eq:normal} is globally asymptotically stable. Consequently, the PID-controlled plant \eqref{eq:plant}--\eqref{eq:pid} satisfies
\begin{equation}
\label{eq:regulation}
    \lim_{t\to\infty}x_1(t)=y^\ast,\qquad
    \lim_{t\to\infty}x_2(t)=0
\end{equation}
for every initial condition and every $f\in F_{L_1,L_2}$.
\end{theorem}

In the original PID parameters, \eqref{eq:new-condition} is
\begin{equation}
\label{eq:new-condition-original}
\begin{aligned}
&\exists\,\mu\in(0,k_d-L_2)\ \text{s.t.}\\
&4\bigl(\mu(k_p-L_1)-k_i\bigr)\\
&\quad >
\mu^2
\Big(
\sqrt{k_d+L_2-\mu}
-\sqrt{k_d-L_2-\mu}
\Big)^2.
\end{aligned}
\end{equation}
Thus, for any $k_i>0$ and $k_d>L_2$, any \(k_p\) for which \eqref{eq:new-condition-original} holds for some \(\mu\) solves the set-point regulation problem robustly.

The same condition can be written as an explicit scalar boundary in the $(p,q)$ plane. For fixed $k>0$, $L\ge 0$, and $q>0$, define
\begin{equation}
\label{eq:Pi-def}
\Pi(q):=
\inf_{0<\mu<q}
\Bigg[
\frac{k}{\mu}
+\frac{\mu}{4}
\Big(\sqrt{q+2L-\mu}-\sqrt{q-\mu}
\Big)^2
\Bigg].
\end{equation}
For each fixed $\mu\in(0,q)$, \eqref{eq:new-condition} is equivalent to requiring that $p$ be larger than the corresponding bracketed expression in \eqref{eq:Pi-def}. Taking the infimum over $\mu$ gives the equivalent boundary form
\begin{equation}
\label{eq:new-boundary}
    p>\Pi(q).
\end{equation}

\begin{remark}
In the special case $L=0$, one obtains
$\Pi(q)=k/q$. For $L>0$, the change of variables
\[
h=\frac{\sqrt{q+2L-\mu}-\sqrt{q-\mu}}
        {\sqrt{q+2L-\mu}+\sqrt{q-\mu}}
\]
transforms the stationarity condition into a fifth-degree polynomial
equation in $h$, so a simple closed-form expression is generally unavailable; nevertheless, the one-dimensional minimization in \eqref{eq:Pi-def} is straightforward numerically.
\end{remark}

For comparison, the scalar sufficient condition of Zhao and Guo \cite{ZhaoGuo2017} is
\begin{equation}
\label{eq:old-zg17-pid}
\begin{aligned}
&(k_p-L_1)(k_d-L_2)\\
&\quad >
k_i+L_2\sqrt{k_i(k_d+L_2)}.
\end{aligned}
\end{equation}
In the notation \eqref{eq:kpq}, this is
\begin{equation}
\label{eq:old-zg17}
    pq>k+L\sqrt{k(q+2L)},
\end{equation}
or equivalently
\begin{equation}
\label{eq:pzg17}
p>p_{17}(q)
:=\frac{k}{q}
+\frac{L}{q}\sqrt{k(q+2L)}.
\end{equation}
Zhang and Guo further provided the following refined sufficient condition for the scalar case \cite{ZhangGuo2019}
\begin{equation}
\label{eq:old-zg19}
    pq>k+L\sqrt{kq},
\end{equation}
or equivalently
\begin{equation}
\label{eq:pzg19}
p>p_{19}(q)
:=\frac{k}{q}
+L\sqrt{\frac{k}{q}}.
\end{equation}

The next proposition establishes the strict ordering of the scalar boundaries.

\begin{proposition}
\label{prop:comparison}
Let $k>0$, $q>0$, and $L>0$. Then
\begin{equation}
\label{eq:comparison-chain}
    \Pi(q)<p_{19}(q)<p_{17}(q).
\end{equation}
If $L=0$, then $\Pi(q)=p_{19}(q)=p_{17}(q)=k/q$.
\end{proposition}

Proposition~\ref{prop:comparison} shows that the region certified by the new sufficient condition strictly contains the scalar sufficient regions obtained in~\cite{ZhaoGuo2017,ZhangGuo2019} whenever $L_2>0$. This strict inclusion is the precise sense in which the present result reduces conservatism. This comparison concerns only the scalar boundaries associated with the sufficient conditions in \cite{ZhaoGuo2017,ZhangGuo2019}; the high-dimensional, multi-agent, and exponential-convergence results of \cite{ZhangGuo2019} remain complementary.

\begin{remark}
\label{rem:no-exponential}
Theorem~\ref{thm:main} establishes global asymptotic stability.
Although the proof below yields a decay estimate for \(y\) and \(z\),
it does not directly yield a decay estimate for the integral state \(x\).
Therefore, no exponential convergence claim is made.
\end{remark}

\section{Proofs of Main Results}
\label{sec:proof}

\subsection{Proof of Theorem~\ref{thm:main}}

Fix $\mu\in(0,q)$ satisfying \eqref{eq:new-condition}, and define
\begin{equation}
\label{eq:C-def}
C:=p+\mu^2
+\mu\sqrt{(q-\mu)(q+2L-\mu)}.
\end{equation}
Consider the Lyapunov candidate
\begin{equation}
\label{eq:V-def}
\begin{aligned}
V(x,y,z)
=&\frac12
\bigl(
\mu kx^2+2kxy+Cy^2\\
&\quad +2\mu yz+z^2
\bigr)
+\int_0^y(\phi(s)-p)s\,ds.
\end{aligned}
\end{equation}

The integral term removes the dependence of the mixed $yz$-coefficient on the unknown function $\phi(y)$, while the coefficient $C$ is chosen to balance the endpoint estimates over the shifted damping interval $r=\psi-\mu$.

\emph{Positive definiteness.}
Since $\phi(s)\ge p$ for all $s$, the integral term in \eqref{eq:V-def} is nonnegative for all $y\in\R$. It remains to check the quadratic part. Its symmetric matrix is
\begin{equation}
\label{eq:P-matrix}
P_\mu=
\begin{pmatrix}
\mu k & k & 0\\
k & C & \mu\\
0 & \mu & 1
\end{pmatrix}.
\end{equation}
Because the right-hand side of \eqref{eq:new-condition} is nonnegative, \eqref{eq:new-condition} implies
\begin{equation}
\label{eq:mupk}
    \mu p-k>0.
\end{equation}
Moreover,
\begin{equation}
\label{eq:Cminus}
C-\mu^2
=p+\mu\sqrt{(q-\mu)(q+2L-\mu)}
>\frac{k}{\mu}.
\end{equation}
The leading principal minors of $P_\mu$ are
\begin{equation}
    \mu k>0,
\end{equation}
\begin{equation}
    \det
    \begin{pmatrix}
        \mu k & k\\
        k & C
    \end{pmatrix}
    =
    k(\mu C-k)>0,
\end{equation}
and
\begin{equation}
\det P_\mu
=k\bigl(\mu(C-\mu^2)-k\bigr)>0.
\end{equation}
Thus $P_\mu$ is positive definite. Hence $V$ is positive definite and radially unbounded.

\emph{Derivative estimate.}
Let $V_0$ denote the quadratic part of $V$. Along \eqref{eq:normal},
\begin{align}
\dot V_0
&=(\mu kx+ky)y+(kx+Cy+\mu z)z \notag\\
&\quad +(\mu y+z)(-kx-\phi y-\psi z) \notag\\
&=(k-\mu\phi)y^2
 +(C-\mu\psi-\phi)yz \notag\\
&\quad +(\mu-\psi)z^2.
\end{align}
The time derivative of the integral term in \eqref{eq:V-def} is
\begin{equation}
\frac{d}{dt}\int_0^y(\phi(s)-p)s\,ds
=(\phi(y)-p)yz.
\end{equation}
Therefore
\begin{equation}
\label{eq:Vdot-exact}
\dot V
=(k-\mu\phi)y^2
+(C-p-\mu\psi)yz
+(\mu-\psi)z^2.
\end{equation}
Using $\phi(y)\ge p$ and $\psi(y,z)\in[q,q+2L]$, we obtain
\begin{equation}
\label{eq:Vdot-bound}
\dot V
\le
-(\mu p-k)y^2
+(C-p-\mu\psi)yz
-(\psi-\mu)z^2.
\end{equation}

It remains to prove that the quadratic form on the right-hand side of \eqref{eq:Vdot-bound} is uniformly negative definite for every $\psi\in[q,q+2L]$. Set
\begin{equation}
r:=\psi-\mu,\quad
a:=q-\mu,\quad
b:=q+2L-\mu.
\end{equation}
Then $0<a\le r\le b$. From \eqref{eq:C-def},
\begin{equation}
    C-p-\mu\psi
    =
    \mu(\sqrt{ab}-r).
\end{equation}
Consequently,
\begin{equation}
\label{eq:ratio}
\frac{(C-p-\mu\psi)^2}{\psi-\mu}
=
\frac{\mu^2(\sqrt{ab}-r)^2}{r}.
\end{equation}
The scalar function
\begin{equation}
    F(r)=\frac{\mu^2(\sqrt{ab}-r)^2}{r}
\end{equation}
satisfies
\begin{equation}
    F'(r)=\mu^2\left(1-\frac{ab}{r^2}\right),
\end{equation}
whose sign changes at $r=\sqrt{ab}$. Hence $F$ is decreasing on $[a,\sqrt{ab}]$ and increasing on $[\sqrt{ab},b]$. Therefore its maximum over $[a,b]$ is attained at the endpoints, and
\begin{equation}
\label{eq:max-ratio}
\begin{aligned}
&\max_{\psi\in[q,q+2L]}
\frac{(C-p-\mu\psi)^2}{\psi-\mu}\\
&\quad =
\mu^2
\Big(
\sqrt{q+2L-\mu}
-\sqrt{q-\mu}
\Big)^2.
\end{aligned}
\end{equation}
By \eqref{eq:new-condition} and \eqref{eq:max-ratio},
\begin{equation}
\label{eq:negdef}
4(\mu p-k)(\psi-\mu)
>(C-p-\mu\psi)^2
\end{equation}
for all $\psi\in[q,q+2L]$. Thus the right-hand side of \eqref{eq:Vdot-bound} is negative definite in $(y,z)$ for each $\psi\in[q,q+2L]$. Equivalently, the negative of this quadratic form is represented by a symmetric matrix whose smallest eigenvalue is positive for each such $\psi$. Since this eigenvalue depends continuously on $\psi$ and the interval $[q,q+2L]$ is compact, it has a positive minimum. Hence there exists $\alpha>0$ such that
\begin{equation}
\label{eq:Vdot-final}
    \dot V\le -\alpha(y^2+z^2)
\end{equation}
for all $(x,y,z)\in\R^3$.

\emph{Convergence.}
Since $V$ is positive definite and radially unbounded, and since $\dot V\le 0$, all maximal solutions are bounded and hence forward complete. By \eqref{eq:Vdot-final}, $\dot V=0$ implies $y=0$ and $z=0$. The largest invariant subset of $\{(x,y,z):y=0,\ z=0\}$ is the origin, because on this set
\begin{equation}
    \dot z=-kx,
\end{equation}
and $k>0$. LaSalle's invariance principle therefore implies global asymptotic stability of the origin of \eqref{eq:normal}. Finally, $y=e=y^\ast-x_1$ and $z=\dot e=-x_2$, so \eqref{eq:regulation} follows.

\subsection{Proof of Proposition~\ref{prop:comparison}}

To avoid confusion with the endpoint-ratio function $F(r)$ used in the proof
of Theorem~\ref{thm:main}, denote the scalar objective in \eqref{eq:Pi-def} by
\begin{equation}
    H(\mu)
    :=
    \frac{k}{\mu}
    +
    \frac{\mu}{4}
    \left(
    \sqrt{q+2L-\mu}-\sqrt{q-\mu}
    \right)^2.
\end{equation}
Then $\Pi(q)=\inf_{0<\mu<q}H(\mu)$. For $0<\mu<q$ and $L>0$,
\begin{equation}
\label{eq:diff-id}
\begin{aligned}
&\sqrt{q+2L-\mu}-\sqrt{q-\mu}\\
&\quad =
\frac{2L}{\sqrt{q+2L-\mu}+\sqrt{q-\mu}}.
\end{aligned}
\end{equation}
Thus
\begin{equation}
\label{eq:HlessG}
H(\mu)
<
\frac{k}{\mu}
+\frac{\mu L^2}{4(q-\mu)}
=:G(\mu).
\end{equation}
Set
\begin{equation}
s:=\sqrt{\frac{k}{q}},\qquad
\mu_0:=q\frac{2s}{L+2s}.
\end{equation}
Then $0<\mu_0<q$, and direct substitution gives
\begin{equation}
\label{eq:Gmu0}
    G(\mu_0)=\frac{k}{q}+L\sqrt{\frac{k}{q}}=p_{19}(q).
\end{equation}
By \eqref{eq:HlessG} and the definition of $\Pi(q)$,
\begin{equation}
\label{eq:Pi-less-p19}
\Pi(q)\le H(\mu_0)<G(\mu_0)=p_{19}(q).
\end{equation}
Moreover,
\begin{equation}
p_{17}(q)-p_{19}(q)=
\frac{L}{q}
\bigl(\sqrt{k(q+2L)}-\sqrt{kq}\bigr)>0,
\end{equation}
which proves \eqref{eq:comparison-chain}. If $L=0$, then \eqref{eq:Pi-def} reduces to $\Pi(q)=\inf_{0<\mu<q} k/\mu=k/q$, and both previous boundaries reduce to $k/q$.

\section{Interpretation}
\label{sec:comparison}

\subsection{Geometric interpretation}

The boundary $p=\Pi(q)$ comes from an interval-wise Lyapunov estimate, not from minimizing a previously known PID inequality. For fixed \(\mu\), set \(a=q-\mu\), \(b=q+2L-\mu\),
and \(r=\psi-\mu\). Before fixing \(C\), the mixed-term estimate
amounts to placing a scalar \(\gamma\), determined by \(C\), in the
normalized worst-case penalty
\[
    \max_{r\in[a,b]}\frac{(\gamma-r)^2}{r}.
\]
Thus the choice of \(C\) is a minimax placement problem over the whole
effective-damping interval. The minimax placement is not the arithmetic
midpoint, but the geometric mean \(\gamma=\sqrt{ab}\), which equalizes
the two endpoint penalties:
\[
    \frac{(\sqrt{ab}-a)^2}{a}
    =
    \frac{(\sqrt{ab}-b)^2}{b}.
\]
With this choice,
\[
\begin{aligned}
C-p-\mu\psi
&= \mu(\sqrt{ab}-r),\\
\frac{(C-p-\mu\psi)^2}{\psi-\mu}
&=
\mu^2\frac{(\sqrt{ab}-r)^2}{r}.
\end{aligned}
\]This explains why the balancing is not merely an algebraic reformulation of the
Lyapunov estimate: it minimizes the uniform normalized mixed-term
penalty induced by the damping interval. This is the precise mechanism
that lowers the scalar sufficient boundary, relative to the certificates
in \cite{ZhaoGuo2017,ZhangGuo2019}, to \(\Pi(q)\).

The same condition can also be viewed in $(p,k,q)$ space. Define
\begin{equation}
\label{eq:Omega-new}
\Omega_{\rm new}:=\{(p,k,q):p>0,\ q>0,\ 0<k<K_{\rm new}(p,q)\},
\end{equation}
where
\begin{equation}
\label{eq:Knew}
\begin{aligned}
K_{\rm new}&(p,q)
\\& := \sup_{0<\mu<q}
\biggl[
\mu p-\frac{\mu^2}{4}\times
\bigl(\sqrt{q+2L-\mu}-\sqrt{q-\mu}\bigr)^2
\biggr].
\end{aligned}
\end{equation}
Then \eqref{eq:new-condition} is equivalent to $(p,k,q)\in\Omega_{\rm new}$. Since the subtracted term in \eqref{eq:Knew} is nonnegative and $\mu<q$, one has $K_{\rm new}(p,q)\le pq$, so the certified sufficient upper surface remains below the necessary upper surface $k=pq$ specified in \eqref{eq:necessary}. Proposition~\ref{prop:comparison} shows that, in the $p$-boundary view, the new condition lowers the scalar sufficient-condition boundaries reported in \cite{ZhaoGuo2017,ZhangGuo2019} from $p_{17}(q)$ and $p_{19}(q)$ to $\Pi(q)$, thereby closing part of the certified gap without claiming optimality for the full nonlinear class.

\begin{figure}[!htb]
\centering
\includegraphics[width=\columnwidth]{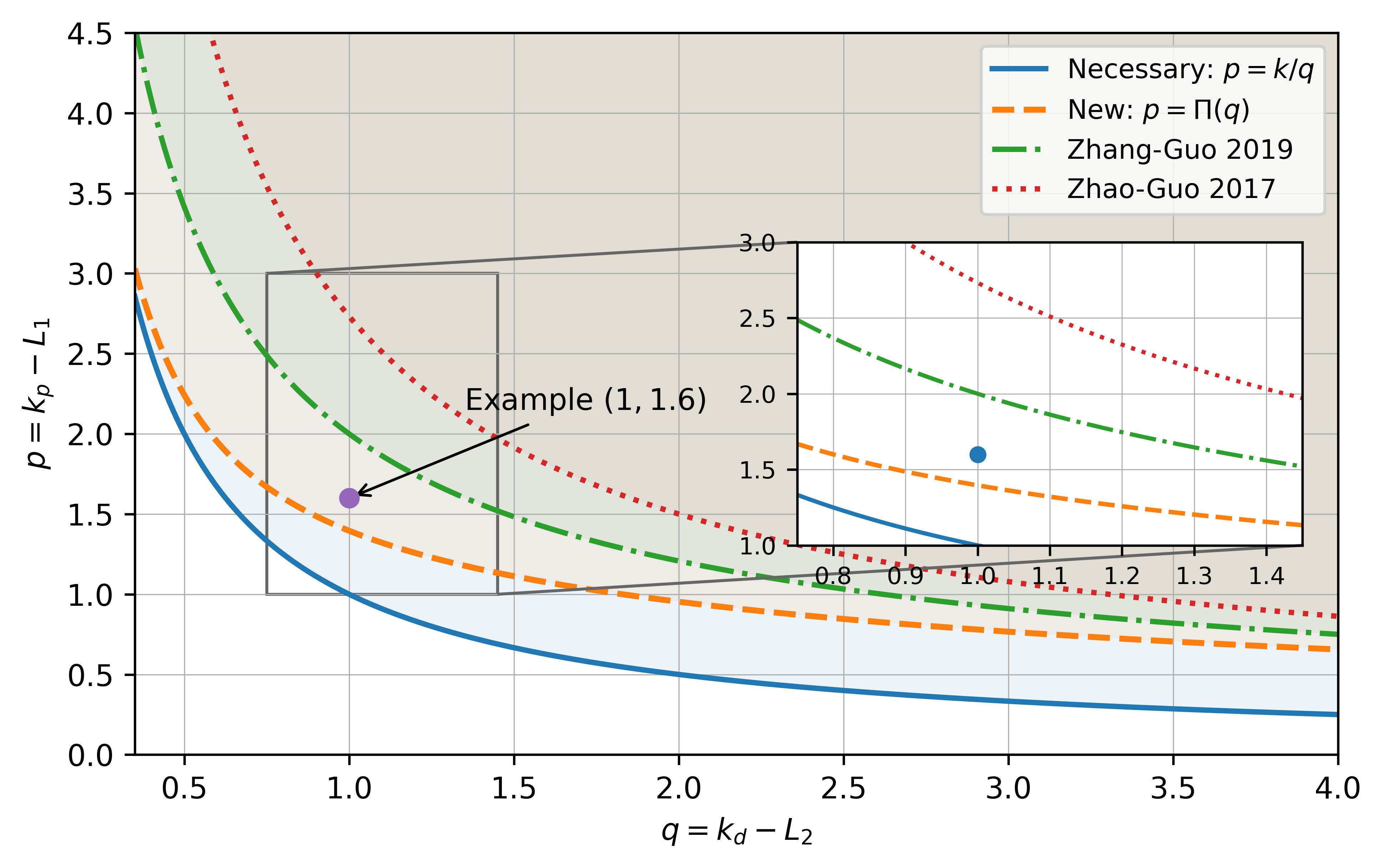}
\caption{Comparison of PID parameter-region boundaries for $k=1$ and $L=1$. The admissible region for each sufficient condition lies above the corresponding boundary. The blue curve is the necessary boundary $p=k/q$, the orange curve is the new sufficient boundary $p=\Pi(q)$, and the green and red curves are the scalar sufficient-condition boundaries in \cite{ZhangGuo2019} and \cite{ZhaoGuo2017}, respectively. The point $(q,p)=(1,1.6)$ satisfies the new condition but violates both previous scalar sufficient conditions.}
\label{fig:parameter-region}
\end{figure}

Figure~\ref{fig:parameter-region} visualizes this comparison for $k=1$ and $L=1$. The regions above the plotted boundaries are admissible for the corresponding sufficient conditions. The point $(q,p)=(1,1.6)$ satisfies the new condition but violates both previous scalar sufficient inequalities.

Table~\ref{tab:boundary-comparison} reports representative boundary values. The last column gives the fraction
$\gamma_{19}(q):=(p_{19}(q)-\Pi(q))/(p_{19}(q)-k/q)\times100\%$
of the certified gap between $p_{19}(q)$ and the necessary boundary $k/q$ that is closed by the new boundary. This is only a certified-boundary metric and should not be interpreted as an optimality claim for the full nonlinear class.

\begin{table}[!htb]
\caption{Numerical comparison of parameter-region boundaries and certified-gap reduction for $k=1$ and $L=1$.}
\label{tab:boundary-comparison}
\centering
\footnotesize
\begin{tabular}{cccccc}
\hline
$q$ & $k/q$ & $\Pi(q)$ & $p_{19}(q)$ & $p_{17}(q)$ & $\gamma_{19}(q)$ \\
\hline
0.5 & 2.000 & 2.242 & 3.414 & 5.162 & 82.9\% \\
1.0 & 1.000 & 1.396 & 2.000 & 2.732 & 60.4\% \\
2.0 & 0.500 & 0.953 & 1.207 & 1.500 & 35.9\% \\
4.0 & 0.250 & 0.656 & 0.750 & 0.862 & 18.9\% \\
\hline
\end{tabular}
\end{table}

\subsection{Relation to the necessary condition}

The worst-case linear model yields the necessary condition
\begin{equation}
\label{eq:necessary}
    p>0,\qquad q>0,\qquad k>0,\qquad pq>k.
\end{equation}
In the original PID parameters, this is
\begin{equation}
\label{eq:necessary-pid}
\begin{gathered}
k_p>L_1,\qquad k_i>0,\\
(k_p-L_1)(k_d-L_2)>k_i.
\end{gathered}
\end{equation}
For an affine-in-velocity subclass, \eqref{eq:necessary} is also sufficient \cite{ZhaoGuo2017}. For the full class $F_{L_1,L_2}$, Theorem~\ref{thm:main} closes part of the certified gap between sufficient and necessary PID parameter regions, but it does not prove that \eqref{eq:necessary} is sufficient.

The remaining gap reflects the freedom of the full nonlinear uncertainty class. In $F_{L_1,L_2}$, the effective damping coefficient $\psi(y,z)$ may move with the state over the whole interval $[q,q+2L]$, without any monotonicity or coupling assumptions beyond the derivative bounds. A Lyapunov certificate based only on these bounds must therefore control the mixed term uniformly over the interval and over all directions $(y,z)$.

This is exactly where the present construction improves on the scalar estimates in \cite{ZhaoGuo2017,ZhangGuo2019}. It does not select a more favorable endpoint or tune the PID gains after a worst-case endpoint has been fixed. Instead, it places the quadratic cross term at the endpoint-balanced position in the shifted damping interval, reducing the uniform mixed-term penalty in \eqref{eq:Vdot-bound}. At the same time, the balancing remains uniform: the coefficient $C$ is state-independent and cannot follow the actual value of $\psi(y,z)$ along trajectories. Letting the Lyapunov cross term adapt to the local damping value would introduce derivatives of $\phi$ or $\psi$ that are not controlled by \eqref{eq:function-class}, which is the main technical obstruction to closing the gap to \eqref{eq:necessary}.

The present result is therefore complementary to \cite{ZhaoWangXue2025}. The impossibility results in \cite{ZhaoWangXue2025} identify certain strongly nonlinear classes that are beyond classical PID stabilization, while the positive results in \cite{ZhaoWangXue2025} use nonlinear PID-type controllers whose gains depend on the state or the regulation error. These results do not provide a less conservative constant-gain PID region for the full derivative-bounded class $F_{L_1,L_2}$. Hence the present result should be read as a refinement of the fixed-gain scalar PID certificates in \cite{ZhaoGuo2017,ZhangGuo2019}, not as a competing nonlinear-PID design.

\section{Numerical Illustration}
\label{sec:numerical}

We use the design point highlighted in Fig.~\ref{fig:parameter-region} and Table~\ref{tab:boundary-comparison}:
\[
L_1=L_2=1,\qquad (k_p,k_i,k_d)=(2.6,1,2).
\]
Then $k=1$, $L=1$, $q=1$, and $p=1.6$. Minimizing \eqref{eq:Pi-def} gives $\mu^\ast\approx0.9276$. For the nearby feasible value $\mu=0.9$,
\[
4(\mu p-k)=1.76
>
1.040\approx
\mu^2\bigl(\sqrt{q+2L-\mu}-\sqrt{q-\mu}\bigr)^2,
\]
so \eqref{eq:new-condition} holds. By contrast, \eqref{eq:old-zg19} requires $pq>2$, and \eqref{eq:old-zg17} requires $pq>1+\sqrt3\approx2.732$, whereas here $pq=1.6$.

For the time-domain illustration, we choose the genuinely nonlinear
velocity-dependent uncertainty
\[
    f(x_1,x_2)=\sin x_1+\log(\cosh x_2).
\]
Then $f_{x_1}=\cos x_1\le1$ and
$|f_{x_2}|=|\tanh x_2|\le1$, so $f\in F_{1,1}$. For
$y^\ast=0$, the normal-form nonlinearity is
$g(y,z)=\sin y-\log(\cosh z)$, and the effective damping
coefficient is
\[
\psi(z)=2+\frac{\log(\cosh z)}{z},\quad z\ne0,\qquad \psi(0)=2.
\]
Thus $\psi(z)$ varies in the admissible interval
$[q,q+2L]=[1,3]$ and is not fixed at a single endpoint.
In particular, the uncertainty is not represented by choosing
$\psi\equiv q$ or $\psi\equiv q+2L$; rather, the effective damping
moves with the state over the admissible interval. This makes the
example consistent with the interval-wise mechanism underlying the
endpoint-balanced certificate.

\begin{figure}[!htb]
\centering
\includegraphics[width=\columnwidth]{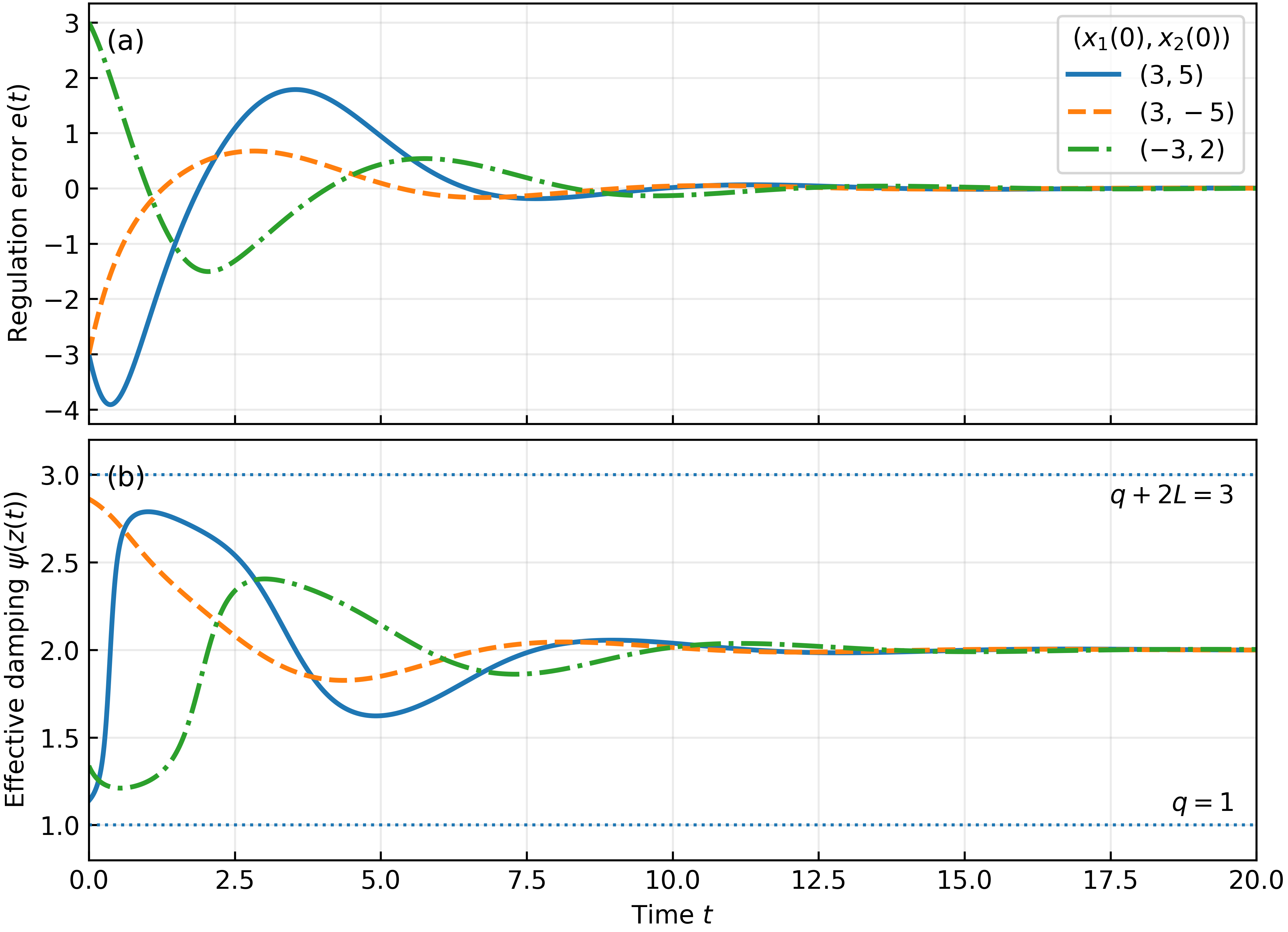}
\caption{Closed-loop behavior for the nonlinear uncertainty
\(f(x_1,x_2)=\sin x_1+\log(\cosh x_2)\) under
\((k_p,k_i,k_d)=(2.6,1,2)\). 
(a) Regulation error \(e(t)=y^\ast-x_1(t)\) with \(y^\ast=0\).
(b) Effective damping coefficient
\(\psi(z(t))=2+\log(\cosh z(t))/z(t)\), with \(\psi(0)=2\) defined by continuity.
The dotted lines show the admissible damping bounds \(q=1\) and
\(q+2L=3\), illustrating damping variation across the interval rather than at a fixed endpoint.}
\label{fig:error-response}
\end{figure}

Figure~\ref{fig:error-response} shows the regulation error and the corresponding effective damping trajectories for initial plant states $(3,5)$, $(3,-5)$, and $(-3,2)$ with $\xi(0)=0$. We integrate and plot the ODEs on $[0,20]$ using RK45 with relative tolerance $10^{-8}$ and absolute tolerance $10^{-10}$. The simulations are consistent with the predicted regulation behavior, while the stability guarantee itself follows from Theorem~\ref{thm:main}.

\section{Conclusion}
\label{sec:conclusion}

This letter has introduced an endpoint-balanced Lyapunov certificate for constant-gain PID control of scalar second-order nonlinear uncertain systems. By balancing the mixed-term penalty at the two endpoints of the admissible effective-damping interval before extracting the scalar PID inequality, the construction yields a strictly larger certified fixed-gain PID region than the regions defined by the scalar sufficient-condition boundaries in~\cite{ZhaoGuo2017,ZhangGuo2019} when \(L_2>0\). For \(L_2=0\), under \(k>0\) and \(q>0\), the corresponding condition recovers the lower-bound form \(p>k/q\) associated with the necessary condition \(pq>k\) obtained from the worst-case linear model. The result establishes global asymptotic set-point regulation for the full uncertainty class defined by the derivative bounds. Future work includes extending the interval-balancing mechanism to high-dimensional uncertain Euler--Lagrange systems and multi-agent settings~\cite{HeLu2023,RoyBaldiIoannou2022,WangZhangBaldiZhong2023,DongChen2022}, and using the certificate as a structured objective in computational searches for Lyapunov functions.

\end{document}